\documentclass[12pt]{amsart}

\usepackage{amsmath}
\usepackage{amscd}
\usepackage{amsfonts}
\usepackage{graphicx}
\usepackage[usenames,dvipsnames]{color}
\usepackage{color}

\newtheorem{theorem}{Theorem}[section]

\newtheorem{lemma}{Lemma}[section]

\newtheorem{definition}{Definition}[section]

\newtheorem{proposition}{Proposition}[section]

\def \epsilon {\varepsilon}
\def \phi {\varphi}
\def \proof {{\bf Proof}\hspace{.4cm}}

\def \Car {Carath\'eodory }

\def \C {\mathbb C}
\def \D {\mathbb D}

\def \cbar {\overline {\mathbb C}} 

\def \dsharpz {|{\mathrm d}^{\scriptscriptstyle \#}\hspace{-.06cm} z|}
\def \i {{\mathbf i}}

\title{Short Separating Geodesics for Multiply Connected Domains}

\author{Mark Comerford}

\address{Department of Mathematics,
University of Rhode Island,
5 Lippitt Road, Room 102F,
Kingston, RI 02881, USA.
email: {\tt mcomerford@math.uri.edu}}



\begin{document}


\renewcommand{\theenumi}{\emph{\arabic{enumi}.}}
\renewcommand{\labelenumi}{\theenumi}

\maketitle

\begin{abstract}
We consider the following questions: given a hyperbolic plane domain and a separation of its complement into two disjoint closed sets each of which contains at least two points, what is the shortest closed hyperbolic geodesic which separates these sets and is it a simple closed curve? We show that a shortest geodesic always exists although in general it may not be simple. However, one can also always find a shortest simple curve and we call such a geodesic a \emph{meridian} of the domain. We prove that, although they are not in general uniquely defined, if one of the sets of the separation of the complement is connected, then they are unique and are also the shortest possible closed curves which separate the complement in this fashion.
\end{abstract}

\section{Introduction}
In this paper we explore the question of how to most efficiently separate the components of the complement of a hyperbolic multiply connected domain using closed geodesics. The study of short closed geodesics, in particular, finding a short homology basis and examining the systole, or shortest geodesic, has been extensively studied in the literature for Riemann surfaces which are compact or at least of finite type \cite{Bus, BS, Par1, Par2}. However, here we will consider exclusively plane domains which, except in rare cases, are not Riemann surfaces of finite type. Also, we will restrict ourselves largely to a fixed homology class. Finally, our motivation for considering these questions is somewhat different as it ultimately comes from certain problems in complex dynamics.

The simple closed geodesics we will obtain, and which we will refer to as \emph{meridians} of our domain, can be used to characterize the conformal type of the domain and give rise to give a short homology basis. They are especially useful in dealing with questions regarding convergence of sequences of pointed domains of the same connectivity with respect to the \Car topology for pointed domains, an example being the question of in which circumstances a limit of pointed 
$n$-connected domains with respect to this topology is again $n$-connected. This in turn is useful for problems in complex dynamics where one needs to be able to uniformize certain fundamental domains in the correct manner, for example in giving the correct generalization of a polynomial-like mapping for a non-autonomous  version of the Sullivan straightening theorem e.g. \cite{Com1}. 

For a hyperbolic domain $U$, denote the hyperbolic metric by $\rho_U (\cdot\,, \cdot)$ or just $\rho (\cdot \,, \cdot)$ if the domain involved is clear from the context. Extending this notation slightly, we shall use $\rho_U(z, A)$ or $\rho(z, A)$ to denote the distance in the hyperbolic metric from a point $z \in U$ to a subset $A$ of $U$. Finally, for a curve $\eta$ in $U$, let us denote the hyperbolic length of $\eta$ in $U$ by $\ell_U(\eta)$, or, again when the context is clear, simply by $\ell(\eta)$.

In this paper, for the sake of convenience, we shall only consider domains which are subsets of $\C$ so that the point at infinity is in one of the components of the complement. However, since all the domains we consider will be hyperbolic and their complements will therefore contain at least three points, there is no loss of generality in doing this. 

We now give a definition which describes three ways in which a curve can be said to divide two sets from each other. For two curves $\gamma_1$, $\gamma_2$ in $U$, we write $\gamma_1 \underset{U}\approx \gamma_2$ to denote homology in $U$.

\begin{definition} Let $E$ and $F$ be two disjoint closed sets and let $\gamma$ be a closed curve which does not meet either set. 

\begin{enumerate}

\item
We say $\gamma$ {\rm separates} $E$ and $F$ if no component of the complement $\cbar \setminus \gamma$ contains points of both $E$ and $F$. 

\vspace{.25cm}
\item
We say $\gamma$ \emph{separates $E$ and $F$ by parity} if the winding number $n(\gamma, z)$ is even for every $z \in E$ and odd for every $z \in F$ or vice versa. 

\vspace{.25cm}
\item
We say $\gamma$ \emph{separates $E$ and $F$ simply} if $\gamma$ separates $E$ and $F$ and at least one of these sets lies in a single component of $\cbar \setminus \gamma$.
\end{enumerate}
\end{definition}

Note that in the first case if $\gamma$ separates $E$ and $F$, then a finite number of components of $\cbar \setminus \gamma$ gives us a separation of $E$ and $F$. 
For a simple closed curve, these three forms of separation are equivalent. If $\gamma$ is a simple closed curve which separates $E$ and $F$, unless stated otherwise, we shall always assume that $\gamma$ has positive orientation and that $n(\gamma, z) = 1$ for all points $z$ in $E$ while $n(\gamma, z) = 0$ for all points $z$ in $F$. 

From now on, for convenience let us assume that, unless stated otherwise, $U$ is a hyperbolic domain and $E$ and $F$ are closed disjoint sets which each contain at least two points and for which $\cbar \setminus U = E \cup F$. Let us call such a separation of the complement of $U$ \emph{non-trivial}. Again for convenience, we shall assume that $E$ is bounded and $\infty \in F$. Note that any curve which separates such sets cannot be homotopic either to a point in $U$ or a point in $\cbar \setminus U$ (i.e. a puncture). 

As we will see later (e.g. in Figure 1), there may be many geodesics in different homotopy classes which separate $E$ and $F$. However, we can always find one which is as short as possible. 
\vspace{.15cm}

\begin{theorem} Let $U$, $E$ and $F$ be as above. Then there exists a geodesic $\gamma$ which separates $E$ and $F$ and whose length in the hyperbolic metric is as short as possible among all geodesics which separate $E$ and $F$. 
\end{theorem}

The two main questions, then, are whether this geodesic is simple and whether it is unique. Unfortunately, the answer to both of these is no as the following result shows. 

\begin{theorem} There exists a hyperbolic domain $U$ and a non-trivial separation $E$, $F$ of $\cbar \setminus U$ so that the shortest geodesic in $U$ which separates $E$ and $F$ is neither simple nor unique.
\end{theorem}

However, on the positive side, we can show the existence of a simple closed geodesic of minimum length which separates the sets $E$, $F$ of a non-trivial separation of the complement of a general hyperbolic domain $U$. 
\vspace{.15cm}

\begin{theorem}
There exists a simple closed geodesic $\gamma$ which separates $E$ and $F$ and whose hyperbolic length is as short as possible for all curves which separate $E$ and $F$ by parity. In particular, $\gamma$ is as short as possible in its homology class and also as short as possible among all simple closed curves which separate $E$ and $F$. Furthermore, any curve which separates $E$ and $F$ by parity and has the same length as $\gamma$ must also be a simple closed geodesic.
\end{theorem}

It is in fact relatively easy to show that there is a shortest closed curve among all simple closed curves which separate $E$ and $F$ and that this curve must be a simple closed geodesic. The strength of this result, then, is that this curve is shortest among a class of curves which are not all simple and that any curve of the same length must also be a simple closed geodesic. 

Let $\gamma$ be a simple closed smooth hyperbolic geodesic which is topologically non-trivial in $U$, let $\pi : \D \mapsto U$ be a universal covering map and let $G$ be the corresponding group of covering transformations.  Any lift of $\gamma$ to $\D$ is a hyperbolic geodesic in $\D$ and going once around $\gamma$ lifts to a hyperbolic M\"obius transformation $A$ which fixes this geodesic. The invariant geodesic is then the axis of $A$, $Ax_A$. The hyperbolic length of $\gamma$ is then the same as the translation length $\ell(A)$ which is the hyperbolic distance $A$ moves points on $Ax_A$. Note that the quantity $\ell(A)$ does not depend on our choice of lift and is conformally invariant. 

We shall call a segment of $Ax_A$ which joins two points $z$, $A(z)$ on $Ax_A$ a \emph{full segment} of $Ax_A$. This discussion and the above result lead to the following definition. 

\begin{definition} Let $U$ be a hyperbolic domain and let $\gamma$ be a simple closed smooth hyperbolic geodesic in $U$ which is as short as possible in its homology class. Such a curve $\gamma$ is called a \emph{meridian} of $U$ and the hyperbolic length $\ell_U(\gamma)$ is called the \emph{translation length} or simply the \emph{length} of $\gamma$. 
\end{definition}

From the statement of Theorem 1.3, a meridian $\gamma$ is the shortest possible simple closed curve which separates $E$ and $F$. An important special case and indeed the prototype for the above definition is the equator of a conformal annulus and just as the equator is important in determining the geometry of a conformal annulus, meridians are important in determining the geometry of domains of (possibly) higher connectivity. In particular, the lengths of the meridians of a domain are conformally invariant. 

The main problem with meridians is that except in special cases such as an annulus, meridians may not be unique. 

\begin{theorem}
There exists a hyperbolic domain $U$ and a non-trivial separation $E$, $F$ of $\cbar \setminus U$ such that there is more than one meridian which separates $E$ and $F$.
\end{theorem}

However, if one of the complementary components is connected, then we do have uniqueness. Furthermore, in this case, the meridian is the shortest possible geodesic which separates $E$ and $F$.
\vspace{.15cm}

\begin{theorem}
If at least one of the sets $E$, $F$ is connected, then there is only one simple closed geodesic $\gamma$ in $U$ which separates $E$ and $F$. In particular, $\gamma$ must be a meridian. In addition, any other geodesic which separates $E$ and $F$ must be longer than $\gamma$.
\end{theorem}

Let us call a meridian as above where at least one of the sets $E$, $F$ is connected a \emph{principal meridian} of $U$. The theorem then tells us that principal meridians are unique. These meridians have other nice properties. For example, we will see that they are disjoint and do not meet any other meridians of $U$ (Theorem 2.5).

\section{Proofs of the Main Results}

In this section, we prove the Theorems stated in Section 1, together with some supporting results.
Before we can prove Theorems 1.1 and especially Theorem 1.3, we need to answer the basic question of whether, given a separation of the complement of a domain into two disjoint closed sets, one can always find at least one simple closed curve which separates these sets
\vspace{.15cm}

\begin{theorem}Let $U$ be a domain and suppose we can find disjoint non-empty closed sets $E$, $F$ with $\cbar \setminus U = E \cup F$. Then there exists a piecewise smooth simple closed curve in $U$ which separates $E$ and $F$.
\end{theorem}

\proof We begin with a variant of the construction given in \cite{Ahl} Theorem 14, Page 139 which is also similar to that described in \cite{New}. As usual, we can suppose that $\infty \in F$ so that $E$ is compact. Let $d < {\rm dist}(E,F)$, the Euclidean distance between $E$ and $F$, and cover the whole plane with a net of squares of side length $\tfrac{d}{\sqrt 2}$. Let $\Gamma$ be the cycle
\vspace{.2cm}
\[\Gamma = \sum_i Q_i\]

where the sum ranges over all those squares $Q_i$ which meet $E$. 
By cancelling those segments which are on the boundary of more than one square, we can say that $n(\Gamma,z) =1$ for every $z \in E$, even if $z$ lies on the edge or corner of a square, while $n(\Gamma, z) = 0$ for every $z \in F$. By examining the possible number of these remaining segments at each vertex, we see that only two or four of them may meet at such a vertex. If we then thicken $\Gamma$ as in \cite{New} in such a way as to avoid $F$, then $\Gamma$ will be non-singular (in the sense defined on Page 140 of \cite{New}) and its boundary will consist of arcs and simple polygons in view of Theorem 1.2 on Chapter 6, Page 140 of \cite{New}. However, the end point of an arc is a vertex which lies on only one segment and so we see that we can express $\Gamma$ in the form
\vspace{.2cm}
\[ \Gamma = \gamma_1 + \gamma_2 + \cdots \cdots + \gamma_n\]

where the curves $\gamma_j$, $1 \le j \le n$ are disjoint simple closed curves. 

We still need to show that we can separate $E$ and $F$ with just one simple closed curve and we do this by induction on the number $n$ of such simple closed curves consisting of horizontal and vertical line segments in a cycle $\Gamma$ as above for which  $n(\Gamma,z) =1$ for every $z \in E$ and $n(\Gamma, z) = 0$ for every $z \in F$. Clearly the result is trivially true for $n = 1$, so suppose now it is true for any such cycle consisting of $n-1$ curves, i.e. that if there is a cycle $\tilde \Gamma$ consisting of 
$n-1$ disjoint simple closed curves which are made up of horizontal and vertical line segments  for which $n(\tilde \Gamma,z) =1$ for every $z \in E$ and $n(\tilde \Gamma, z) = 0$ for every $z \in F$, then we can find a simple closed curve in $U$ consisting of horizontal and vertical line segments and which separates $E$ and $F$.

So suppose now that $\Gamma$ consists of $n$ curves and  
join $\gamma_1$ to $\gamma_2$ by a path $\eta$ in $U$. It is easy to see that we may assume that $\eta$ is an arc (i.e. a non self-intersecting path) consisting of horizontal and vertical line segments. Such a path may cross both these curves a number of times and may also meet some of the other curves in our cycle. However, by restricting the path and relabelling the curves if necessary, we may assume that $\eta$ is a path from $\gamma_1$ to $\gamma_2$ which meets $\gamma_1$ and $\gamma_2$ only at its endpoints and avoids the remaining curves $\gamma_3, \ldots, \gamma_n$. If we then let $\Gamma'$ be the new cycle $\Gamma + \eta - \eta$, then we still have that $n(\Gamma',z) =1$ for every $z \in E$ and $n(\Gamma', z) = 0$ for every $z \in F$. 

By turning $\eta$ and $-\eta$ into two slightly different disjoint arcs in opposite directions (again consisting of horizontal and vertical line segments) and deleting a small part of $\gamma_1$ and $\gamma_2$, we can then combine $\gamma_1$ and $\gamma_2$ into one simple closed curve in $U$ which is disjoint from $\gamma_3, \ldots, \gamma_n$. 
(Note that we can make these two paths by replacing $\eta$ with a suitable collection of small squares whose boundaries are oriented positively and then cancelling those segments which lie on the boundary of more than one square including those small segments of $\gamma_1$ and $\gamma_2$ that we remove above).

If we denote the new cycle obtained in this way by $\Gamma''$, then we can ensure that this cycle is homotopic and thus homologous in $U$ to $\Gamma'$ and so we still have that $n(\Gamma'',z) =1$ for every $z \in E$ and $n(\Gamma'', z) = 0$ for every $z \in F$. $\Gamma''$ then consists of $n-1$ simple closed curves and it follows by induction that we can find a single simple closed curve $\gamma$ in $U$ consisting of horizontal and vertical line segments and which separates $E$ and $F$ as required. $\Box$

We remark that using this result, one can easily satisfy oneself that finitely connected domains are precisely those domains where any subset of the components of the complement can be separated from the remaining components by a simple closed curve. For more on finitely connected domains, we refer the reader to Section 3 of this paper.

For the most part we shall be working with the spherical metric $\rm d^{\scriptscriptstyle \#} (\cdot \,, \cdot)$ on $\cbar$ (rather than the Euclidean metric). Recall that the length element for this metric, $\dsharpz$ is given by
\vspace{.25cm}
\[ \dsharpz = \frac{|{\rm d}z|}{1 + |z|^2}.\]

We recall a well-known lemma concerning the behaviour of the hyperbolic metric near the boundary. A proof can be found in \cite{CG} Page 13, Theorem 4.3 and it is the lower bound it gives on the hyperbolic metric which will be of particular importance for us. For a point $u \in U$, we shall denote the spherical distance to $\partial U$ by $\delta^\#_U (u)$ or just $\delta^\# (u)$ if once again the domain is clear from the context. For a family $\Phi = \{\phi_\alpha, \alpha \in A\}$ of M\"obius transformations, we say that $\Phi$ is \emph{bi-equicontinuous on $\cbar$} if both $\Phi$ and the family $\Phi^{\circ -1} =  \{\phi_\alpha^{\circ -1}, \alpha \in A\}$ of inverse mappings are equicontinous families on $\cbar$ (with respect to the spherical metric). 

\begin{lemma} Let $U \subset \C$ be a hyperbolic domain. Then there exists $K \ge 1$ for which the hyperbolic metric $\rho(\cdot\,, \cdot)$ on $U$ satisfies
\[  \frac{1}{K}\frac{1 + o(1)}{\delta^\#(z)\log(1/\delta^\#(z))} \dsharpz \le {\rm d}\rho (z) \le \frac{K}{\delta^\#(z)} \dsharpz, \qquad \mbox{as} \quad z \to \partial U.\]
\end{lemma}

An important fact to note in connection with this result is the degree to which these estimates are independent of the particular domain under consideration and in particular how we can get bounds on the constant $K$ and even more importantly on the size of the region on which this estimate holds near a given point on the boundary. Let $z_1$ be the closest point in $\cbar \setminus U$ to $z$ and choose two other points $z_2$, $z_3$ in $\partial U$. whose separation in the spherical topology is bounded below. The proof in \cite{CG} then uses a M\"obius transformation to map these three points to $0$, $1$ and $\infty$ respectively and uses the Schwarz lemma to compare the hyperbolic metric for $U$ with that of $\cbar \setminus \{0,1,\infty\}$. It thus follows that these estimates are uniform with respect to the minimum separation in the spherical metric between $z_1$, $z_2$ and $z_3$. This fact will be important to us later on as we will see in the counterexample in Theorem 1.2. Finally, note that the original statement of this result is in terms of  the Euclidean rather than the spherical distance to the boundary. However, in view of the above discussion, it is clear that if we want our estimates to be uniform, we must make use of the spherical distance.

We remind the reader of two well-known results on hyperbolic geodesics.

\begin{theorem}[\cite{Bus} Theorem 1.5.3 i), vi)] Let $U$ be a hyperbolic domain and let $z_1, z_2$ be two (not necessarily distinct) points of $U$. Then there is a unique geodesic in every homotopy class of curves joining $z_1$ and $z_2$ which is also the unique shortest curve in this homotopy class.
\end{theorem}

\begin{theorem}[\cite{KL} Theorem 7.2.5] Let $U$ be a hyperbolic domain and let $\gamma$ be a non-trivial closed curve in $U$. Then either there is a unique shortest geodesic which is also the shortest curve in the free homotopy class of $\gamma$, or $\gamma$ is homotopic to a puncture on $U$ and there are arbitrarily short curves in the free homotopy class of $\gamma$. 
\end{theorem}

We first consider the easier question of the shortest curve in a given homotopy class which separates $E$ and $F$.

\begin{theorem} Let $\tilde \gamma$ be a simple closed curve which separates $E$ and $F$. Then there exists a unique simple closed smooth geodesic $\gamma$ which is the shortest curve in the free homotopy class of $\tilde \gamma$ in $U$ and in particular also separates $E$ and $F$. 

Conversely, given a simple closed smooth hyperbolic geodesic $\gamma$ in $U$, $\gamma$ separates $\cbar \setminus U$ non-trivially and is the unique geodesic in its homotopy class and also the unique curve of shortest possible length in this class. 
\end{theorem}

Note that the fact that $\gamma$ must separate $E$ and $F$ in the first part of the statement follows easily from the Jordan curve theorem and the fact that $\gamma$ is simple and must be homologous in $U$ to $\tilde \gamma$.

{\bf Proof of Theorem 2.4 \hspace{.4cm}} We start by observing that in view of Theorem 2.1, we can always find a simple closed curve $\tilde \gamma$ which separates $E$ and $F$. Since $E$ and $F$ both contain at least two points, any curve which is homotopic to $\tilde \gamma$ cannot be homotopic to a point or a puncture on $U$ and so by Proposition 3.3.9 on Page 73 of \cite{Hub}, then there is a simple closed geodesic $\gamma$ in the free homotopy class of $\tilde \gamma$ which, using winding numbers, must separate $E$ and $F$. By Theorem 2.3, $\gamma$ is unique. 

To prove the second part of the statement, let $\gamma$ be a simple closed hyperbolic geodesic in $U$ and as above let $E$ and $F$ denote the two subsets of $\cbar \setminus U$ separated by $\gamma$. If one of these sets, say $E$, is empty, then it is a routine exercise to show that $\gamma$ is contractible to a point in $U$ and so cannot be a geodesic. Also, if one of these sets, again say $E$, is a single point, then once more it is easy to see that $\gamma$ is homotopic in $U$ to a puncture and so again by Theorem 2.3 cannot be a geodesic which completes the proof. $\Box$

A more general question is to consider all geodesics which separate $E$ and $F$ as in the first part of Definition 1.1 and this is the content of Theorem 1.1.

{\bf Proof of Theorem 1.1 \hspace{.4cm}} We note that, by Theorems 2.1 and 2.4, the class of geodesics which separate $E$ and $F$ is non-empty. Since neither $E$ nor $F$ is a point, the spherical and hence Euclidean diameter and hence also the arc length of any curve which separates these sets must be bounded below. 

Now let $l$ be the infimum of the hyperbolic lengths (in the sense of the Riemann-Stieltjes integral e.g. in Chapter IV of \cite{Con}) of all closed geodesics which separate $E$ and $F$ and take a sequence $\gamma_n$ of such geodesics whose hyperbolic lengths tend to $l$ (note that these curves are smooth and thus of bounded variation by Proposition 1.3 on Page 58 of \cite{Con}). Since the Euclidean arc lengths of these curves are bounded below, the standard estimates on the hyperbolic metric from Lemma 2.1 combined with the remarks following it and the divergence of the improper integral
\vspace{.16cm}
\begin{equation} \int_0^{\tfrac{1}{2}} { \frac{1}{x \log (1/ x)}{\rm d}x}\end{equation}

show that these curves are uniformly bounded away from $E$ and $F$. The curves $\gamma_n$ then all lie in a compact subset of $U$ and as their hyperbolic lengths are uniformly bounded above, their arc lengths will also be uniformly bounded above. Finally, since their arc lengths are uniformly bounded below, it follows that their hyperbolic lengths are uniformly bounded below away from zero and so $l >0$.

Since the arc lengths of the curves $\gamma_n$ are then uniformly bounded above and below we can use linear maps to reparametrize these curves so that are all defined on $[0,1]$ and give us an equicontinuous family on this interval. Using the Arzela-Ascoli theorem and the completeness of ${\mathbf C}[0,1]$, we can find a subsequence $n_k$ such that as $k$ tends to infinity, the curves $\gamma_{n_k}$ converge uniformly to a limit $\gamma$. 

It follows from the uniform convergence and that fact that the curves $\gamma_{n_k}$ have uniformly bounded variation that $\gamma$ must also be of bounded variation. 
Using equicontinuity and again the fact that the curves $\gamma_{n_k}$ have uniformly bounded variation and are uniformly bounded away from $\partial U$, a standard argument similar to that in the proof of Theorem 1.4 on Page 60 of \cite{Con} shows that 
$\ell(\gamma) = l$. Additionally, as $\gamma$ is a uniform limit of geodesics, it must also be a geodesic in the hyperbolic metric of $U$.

The uniform convergence of $\gamma_{n_k}$ to $\gamma$ ensures that $\gamma_{n_k}$ is homotopic in $U$ to $\gamma$ for $k$ large enough. It then follows immediately from Theorem 2.3 that $\gamma_{n_k} = \gamma$ for $k$ large enough and in particular that $\gamma$ must separate $E$ and $F$ as desired. $\Box$ 

A stronger version of this result would be if $\gamma$ were the shortest curve among all curves which separate $E$ and $F$ and not just the shortest geodesic. A similar procedure to that above shows that it is not difficult to deduce the existence of a closed curve which separates $E$ and $F$ and which is as short as possible. The difficulty is that this curve may not be a geodesic and if one applies Theorem 2.3 to find the shortest geodesic in its free homotopy class, then there is nothing to prevent this geodesic from no longer separating $E$ and $F$.

We next prove Theorem 1.2. First, however, we need to prove a lemma about the points of intersection of hyperbolic geodesics.

\begin{lemma} Let $U$ be a hyperbolic domain and let $\gamma$ be a closed hyperbolic geodesic in $U$. Then $\gamma$ can have only finitely many points of self-intersection. This includes the case of a  point on $\gamma$ which meets infinitely many distinct arcs of $\gamma$.

Further, if $\gamma_1$ and $\gamma_2$ are two distinct closed geodesics in $U$, then $\gamma_1$ and $\gamma_2$ can have only finitely many points of intersection. 
\end{lemma}

\proof Any geodesic in $U$ is easily seen to be freely homotopic in $U$ to a curve consisting of finitely many horizontal and vertical line segments. As two such polygonal curves can intersect each other only finitely many times, the result for two curves follows easily from Proposition 3.3.9 on Page 73 of \cite{Hub}. A similar argument to the proof of this result gives us the case with one geodesic. $\Box$

{\bf Proof of Theorem 1.2 \hspace{.4cm}} The construction of the required counterexample is based on the figure below. 

\vspace{.2cm}
\scalebox{0.378}{\includegraphics{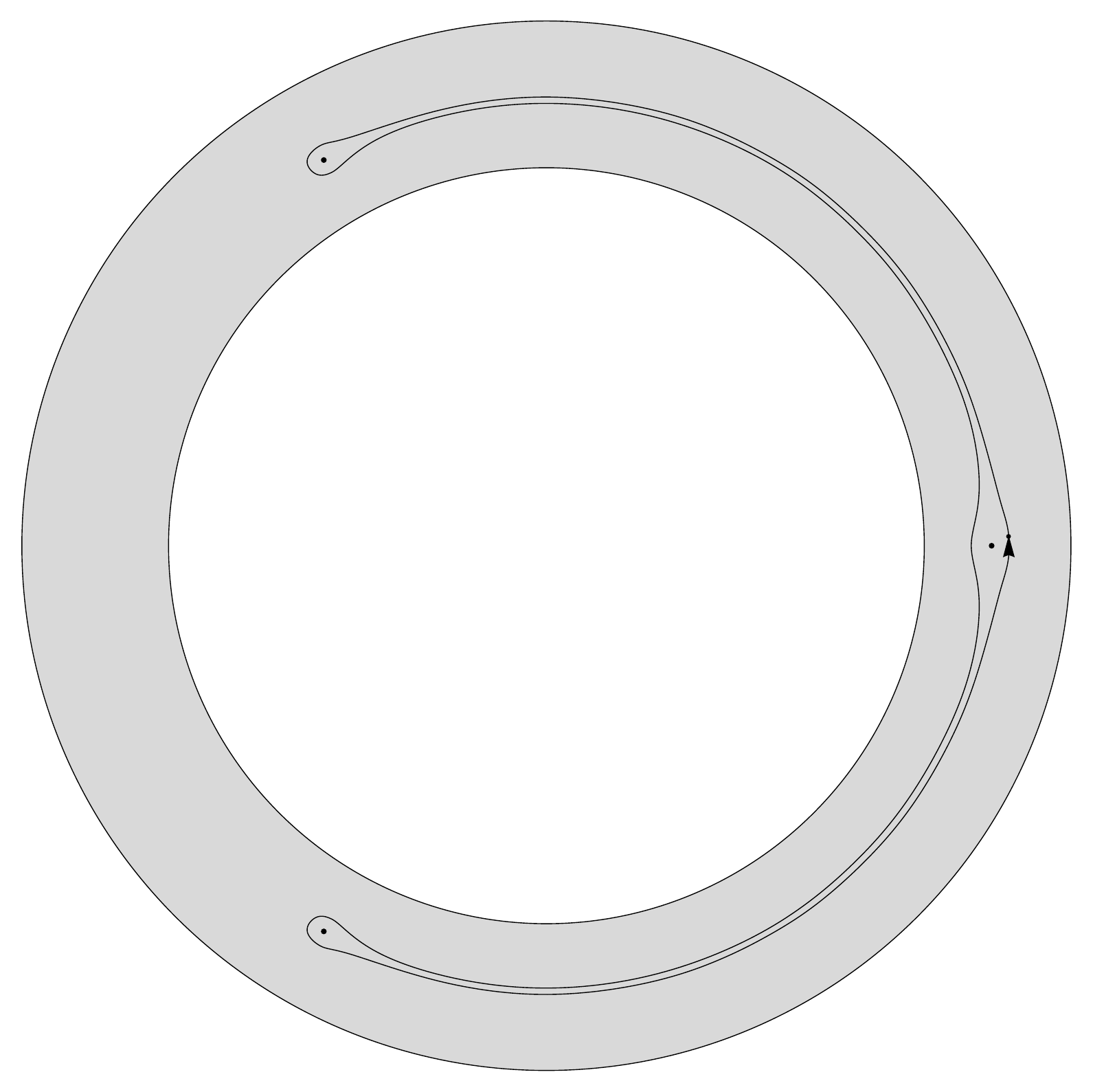}, \hspace{3.5cm}\includegraphics{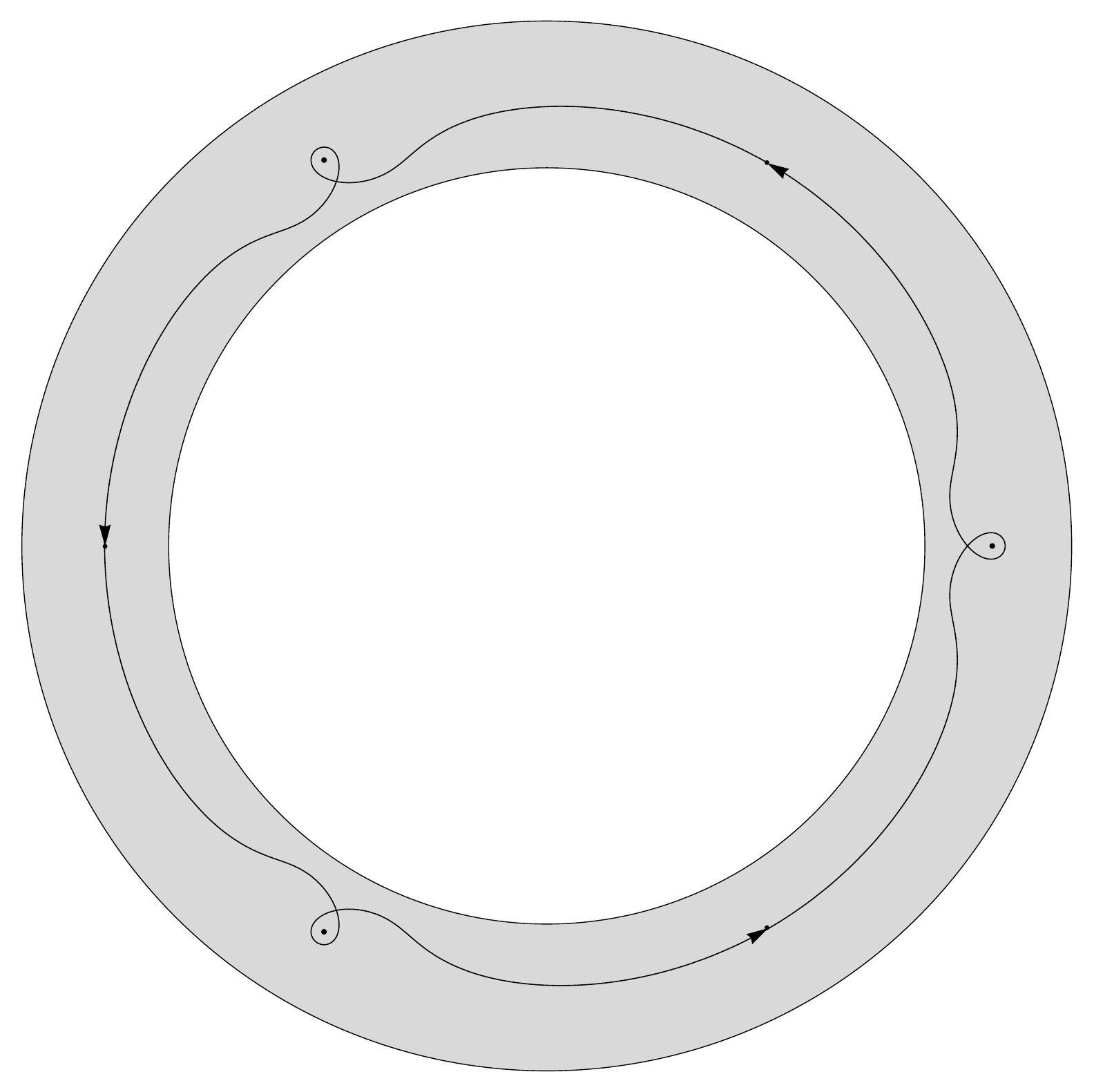}}

\vspace{.4cm}

Let $0 < \epsilon < 1$, let $z^k_\epsilon$ be the three points $(1 + \epsilon) e^{2\pi i k/3}$, $k = 0,1,2$, let $E_\epsilon = \{z^0_\epsilon, z^1_\epsilon,z^2_\epsilon \}$, $F_\epsilon = \overline \D \cup (\cbar \setminus {\mathrm D}(0, (1 + \epsilon)^2))$ and let $U_\epsilon$ be the domain $\cbar \setminus (E_\epsilon \cup F_\epsilon)$. By Lemma 2.1, if we let $\gamma^1_\epsilon$ be a simple closed geodesics which separates $E_\epsilon$ from $F_\epsilon$, then $\gamma^1_\epsilon$ must have hyperbolic length $l^1_\epsilon$ which tends to infinity as $\epsilon$ tends to $0$ (in fact, $l^1_\epsilon = (C + o(1))/\epsilon$ as $\epsilon \to 0$ for some $C > 0$).

We claim that we can find $c >0$ not depending on $\epsilon$ (or the homotopy class of $\gamma^1_\epsilon$) such that the spherical distance from $\gamma^1_\epsilon$ to $\partial U_\epsilon$ is bounded below by $c \epsilon$. To see this note that if $\sigma_\epsilon$ is a segment of $\gamma_\epsilon$ which subtends an angle $\epsilon$ at $0$ and every point on $\sigma_\epsilon$ is within distance $\tfrac{\epsilon}{n}$ of $\partial U_\epsilon$, then by Lemma 2.1 and the remarks following it, we can find $C_1 >0$ independent of $\epsilon$ such that $\ell_{U_{\epsilon}}(\sigma_\epsilon) \ge \tfrac{C_1 n}{\log n}$. On the other hand, if we replace $\sigma_\epsilon$ by a curve $\tau_\epsilon$ consisting of two radial line segments from the endpoints of $\sigma_\epsilon$ to ${\mathrm C}(0, 1 + \epsilon)$ combined with a suitable segment of this circle which subtends angle $\epsilon$, then Lemma 2.1 again tells us that we can find $C_2 >0$ not depending on $\epsilon$ such that $\ell(\tau_\epsilon) \le C_2 \log n$ (note that such a modified curve is still homotopic to $\gamma^1_\epsilon$).

It then follows that we can find $c > 0$ not depending on $\epsilon$ such that for any  segment $\sigma_\epsilon$ which subtends an angle $\ge \epsilon$, we must have a point $z$ on $\sigma_\epsilon$ with $\delta^\#_{U_\epsilon}(z) \ge c \epsilon$. Further, for $\epsilon < 1$, any such point is within distance $\le 7 \epsilon$ of another such point $w$ on $\gamma^1_\epsilon$ with $\delta^\#_{U_\epsilon}(w) \ge c \epsilon$ and we can find an arc of $\gamma^1_\epsilon$ joining $z$ and $w$ which subtends an angle of at most $2\epsilon$. By Lemma 2.1 once more, the hyperbolic distance $\rho_{U_\epsilon} (z, w)$ is uniformly bounded in $\epsilon$ and, as it subtends an angle of at most $2\epsilon$, a simple homotopy argument using Theorem 2.2 shows that the same must be true for the above arc of $\gamma^1_\epsilon$. A final application of Lemma 2.1 together with the divergence of the improper integral in the last result shows that we can make $c$ smaller if needed, but still independent of $\epsilon$ so that $\delta^\#_{U_\epsilon}(z) \ge c\epsilon$ for every point $z$ on $\gamma^1_\epsilon$ as desired. 

Any radial half-line from $0$ which meets the interior of $\gamma^1_\epsilon$ must meet $\gamma^1_\epsilon$ itself in at least two places. Then, by connecting the points of $E_\epsilon$ by curves inside $\gamma^1_\epsilon$, and using the intermediate value theorem, we see that we can find at least four disjoint arcs of $\gamma^1_\epsilon$, each of which covers angles precisely in the range $[0,\tfrac{2\pi}{3}]$, $[\tfrac{2\pi}{3}, \tfrac{4\pi}{3}]$, or $[\tfrac{4\pi}{3}, 2\pi]$.

We can thus choose an arc $\eta_\epsilon$ which subtends one of these intervals and ensure that $\ell(\eta_\epsilon) \le \tfrac{\ell(\gamma^1_\epsilon)}{4}$. From above, and using Lemma 2.1, we can if needed add a radial line segment and part of a circle whose centre lies on ${\mathrm C}(0, 1+ \epsilon)$ to one of the endpoints to obtain a new curve $\tilde \eta_\epsilon$ which is still an arc, subtends the same range of angles as $\eta_\epsilon$ and whose endpoints lie on the same circle about $0$ and such that we can find $C_3 >0$ independent of $\epsilon$ such that $\ell (\tilde \eta_\epsilon) \le \ell (\eta_\epsilon) + C_3$.

We can use $\tilde \eta_\epsilon$ and its rotations by $e^{2\pi i /3}$, $e^{4\pi i /3}$ to give us a simple closed curve $\alpha_\epsilon$ which by symmetry encloses either $\overline \D$ or $E_\epsilon \cup \overline \D$. In the latter case, we may apply the transformation $z \mapsto \tfrac{(1 + \epsilon)^2}{z}$ which preserves hyperbolic length in $U_\epsilon$ so as to assume without loss of generality that $\alpha_\epsilon$ encloses only $\overline \D$ and, by reversing the direction if needed, we can assume that it does so with positive orientation, i.e. the winding number about points of $\overline \D$ is $1$. 
 
Let $C^i_\epsilon$ be circles of radius $\epsilon/2$ about the three points $z^i_\epsilon$, $i = 1,2,3$ which we will orient \emph{negatively} (i.e. clockwise).
We can then combine $\alpha_\epsilon$ with these three circles and, if needed, three short line segments (traversed in both directions) to make a closed curve $\beta_\epsilon$ which separates $E_\epsilon$ and $F_\epsilon$. As $\gamma^1_\epsilon$ is at least distance $c\epsilon$ away from $\partial U_\epsilon$, once more using Lemma 2.1, the contributions of the circles and the line segments are bounded in $\epsilon$ while the contributions of the segments of $\gamma^1_\epsilon$ have length at most $l^1_\epsilon /4$. Thus the total length of this curve is at most 
$C_4 + (3/4)(l^1_\epsilon)$ for some constant $C_4 > 0$ which does not depend on $\epsilon$. If we then apply Theorem 2.3 and let $\gamma^2_\epsilon$ be the shortest closed geodesic in the homotopy class of $\beta_\epsilon$, then $n(\gamma^2_\epsilon, z) = -1$ for the three points of $E_\epsilon$. Also, for $F_\epsilon$, $n(\gamma^2_\epsilon, z) = 1$ for the bounded component $\overline \D$ and $n(\gamma^2_\epsilon, z) = 0$ for the unbounded component $\cbar \setminus {\mathrm D}(0, (1 + \epsilon)^2)$.

By the Jordan curve theorem, $\gamma^2_\epsilon$ cannot then be a simple curve and since its length is at most $C_4 + (3/4)(l^1_\epsilon)$, we see on making $\epsilon$ small enough that $\gamma^2_\epsilon$ is shorter than $\gamma^1_\epsilon$. Thus a simple closed geodesic (and in particular a meridian) is not the shortest curve to separate $E_\epsilon$ and $F_\epsilon$ in this case. 

This proves the first part of the statement. Now let $\gamma_\epsilon$ be the shortest closed geodesic which separates $E_\epsilon$ and $F_\epsilon$. Suppose for the sake of contradiction that $\gamma_\epsilon$ is unique. Since $U_\epsilon$ is symmetric under the transformation, $z \mapsto \tfrac{(1 + \epsilon)^2}{z}$, $\gamma_\epsilon$ must also possess this symmetry. Now, by Lemma 2.2, $\gamma_\epsilon$ can meet ${\mathrm C}(0, 1 + \epsilon)$ at only finitely many points since by symmetry any such point of intersection would necessarily also be a point of self-intersection of $\gamma_\epsilon$ itself. Thus, with only finitely many exceptions, any radial half-line from $0$ which meets $\gamma_\epsilon$ must do so at at least two points. 

On the other hand, $\gamma_\epsilon$ must subtend an angle of at least $\tfrac{4\pi}{3}$ at $0$ since otherwise it could not separate $E_\epsilon$ and $F_\epsilon$. From above,  we can then find at least four arcs as before which meet at only finitely many points and which subtend angles in the ranges $[0,\tfrac{2\pi}{3}]$, $[\tfrac{2\pi}{3}, \tfrac{4\pi}{3}]$, or $[\tfrac{4\pi}{3}, 2\pi]$ and pick one of them, $\eta_\epsilon$, such that $\ell(\eta_\epsilon) \le \tfrac{\ell(\gamma_\epsilon)}{4}$. The same argument as above then shows that, for $\epsilon$ small enough, we can obtain a shorter curve which still separates these sets and, with this contradiction, we see that $\gamma_\epsilon$ cannot be unique as desired. 
$\Box$

{\bf Proof of Theorem 1.3 \hspace{.4cm}} By Theorem 2.1, the class of curves which separate $E$ and $F$ by parity is non-empty. Note that since $E$ and $F$ both contain at least two points, any curve which separates these sets by parity has spherical diameter and hence Euclidean diameter and arc length which are bounded below. Note also that since we have assumed $\infty \in F$ and the winding number of any curve about $\infty$ is zero, for 
any curve which separates $E$ and $F$ by parity, the winding number about points of $F$ will always be even. Now let $l$ be the infimum of the hyperbolic lengths (again in the sense of Riemann Stieltjes) of all closed curves which separate $E$ and $F$ by parity. The same argument as in the proof of Theorem 1.1 shows that $l > 0$. 
 
We then similarly to before take a sequence of curves $\gamma_n$ which separate $E$ and $F$ by parity and whose hyperbolic lengths tend to $l$. By Lemma 2.1 and the same estimate (1) as used in the proof of Theorem 1.1, these curves are clearly uniformly bounded away from $E$ and $F$ and also of uniformly bounded variation (i.e. arc length). 
We can thus again parametrize these curves using arc length and 
use linear maps to reparametrize them so that they are all defined on $[0,1]$ and give us an equicontinuous family on this interval. 

Again by the Arzela-Ascoli theorem and the completeness of ${\mathbf C}[0,1]$, we can find a subsequence $n_k$ such that as $k$ tends to infinity, the curves $\gamma_{n_k}$ converge uniformly to a limit $\gamma$. $\gamma$ must also then be of bounded arc length and the same argument as before shows we must have $\ell(\gamma) = l$. 

We need to check are that $\gamma$ is a simple closed geodesic which separates $E$ and $F$ by parity. However, the uniform convergence of the curves $\gamma_{n_k}$ to $\gamma$ implies that for any $z \in \cbar \setminus U$, the winding numbers $n(\gamma_{n_k}, z)$ must eventually be the same as $n(\gamma, z)$ and with this it is obvious that $\gamma$ does indeed separate these sets by parity as desired and since $E$ and $F$ both contain at least two points, $\gamma$ cannot be homotopic to a point or a puncture (since otherwise the winding number of $\gamma$ will be the same about at least three points of $\cbar \setminus U$ where at least one of these points comes from $E$ and at least one from $F$). Any curve in the homotopy class of $\gamma$ will also separate $E$ and $F$ by parity and so, using the minimality of the length of $\gamma$, by the uniqueness part of Theorem 2.3, $\gamma$ must be a geodesic.

So suppose now that $\gamma$ is not simple and thus intersects itself. If $w$ is a point of self-intersection, then $\gamma$ can be expressed as the composition of two different closed curves $\gamma_1$ and $\gamma_2$ both beginning and ending at $w$.

Let $\tilde \gamma$ be the curve $\gamma_2 \circ -\gamma_1$. Note that $\tilde \gamma$ has the same length as $\gamma$ and since there can only be one geodesic in any direction through a given point, any self-intersections of $\gamma$ will be transversal and so $\tilde \gamma$ cannot be a geodesic as it must have (at least two) corners (i.e. different tangents to the curve at the same point when approached from different directions). Now if $z \in \cbar \setminus U$, $n(\tilde \gamma, z) = n(\gamma, z) - 2n(\gamma_1, z)$ from which it follows that $\tilde \gamma$ also separates $E$ and $F$ by parity and thus is not homotopic to a point or a puncture. If we now apply Theorem 2.3 again to obtain a geodesic $\tilde {\tilde \gamma}$ which is homotopic to $\tilde \gamma$, then, as $\tilde \gamma$ is not a geodesic, $\tilde {\tilde \gamma}$ is shorter than $\tilde \gamma$ and thus than $\gamma$. With this contradiction, it follows that $\gamma$ is simple as desired. 

Finally, if $\gamma'$ is any other curve of length $l$ which separates $E$ and $F$ by parity, any curve in the free homotopy class of $\gamma'$ must also separate $E$ and $F$ by parity and it follows from Theorem 2.2 and the minimality of $l$ that $\gamma'$ must be a geodesic and from above that $\gamma'$ must then be simple, which completes the proof. $\Box$ 

As we observed earlier, the geodesic $\gamma$ above is as short as possible among a set of curves which may be self-intersecting. If one only considers simple closed curves, then it is possible to use Theorem 2.4 to give an easier proof using a similar argument to that for Theorem 1.1. 

It is also worth noting that, in the proof of Theorem 1.2, the curve $\gamma^2_\epsilon$ does not separate $E_\epsilon$ and $F_\epsilon$ by parity. Nor does it separate these sets simply.

{\bf Proof of Theorem 1.4 \hspace{.4cm}}The proof is based on the following figure.

\scalebox{0.82}{\includegraphics{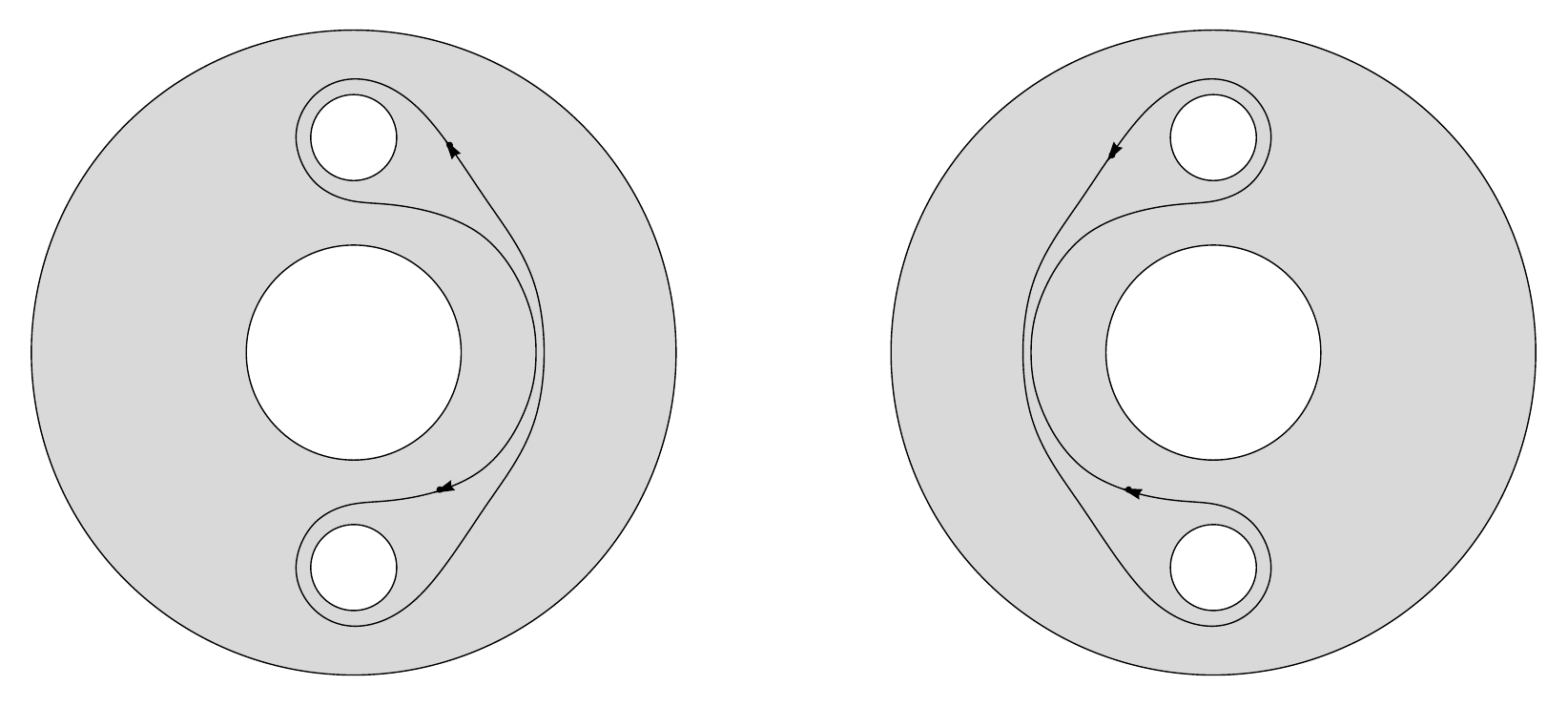}}

Let $E = \overline {{\mathrm D}(3.5 \i, .5)} \cup \overline {{\mathrm D}(-3.5 \i, .5)}$, $F = \overline{{\mathrm D}(0, 2)} \cup (\cbar \setminus {\mathrm D}(0,5))$ and set $U = \cbar \setminus (E \cup F)$. Let $\gamma$ be a meridian which separates $E$ and $F$. The quadruply connected domain $U$ is then symmetric under $T: z \mapsto -z$ and if $\gamma$ were unique, it would also possess this symmetry. 

Let $w$ be a point on $\gamma$ and let $\sigma$ be an arc of $\gamma$ which connects $w$ to $-w$. If we then set $\tilde \gamma = T(\sigma) \circ \sigma$, then $\tilde \gamma$ is a closed curve. The change in argument about $0$ by integrating $\tfrac{1}{z}$ along $\sigma$ and along $T(\sigma)$ must be the same odd multiple of $\pi$. 
It follows that the winding number of $\tilde \gamma$ about points of $\overline{{\mathrm D}(0, 2)}$ must be non-zero. Since the track of $\tilde \gamma$ is a subset of that of $\gamma$, $\gamma$ must then
separate $\overline{{\mathrm D}(0, 2)}$ from $\cbar \setminus {\mathrm D}(0,5)$ which is impossible. $\Box$

We have already seen that the meridian is not in general the shortest curve which separates $E$ and $F$. We next examine whether it is the shortest curve which separates $E$ and $F$ simply.  

\begin{theorem} Any geodesic of minimum length which separates $E$ and $F$ simply must be a meridian. 
\end{theorem}

\proof As at the start of the proofs of Theorems 1.1 and 1.2, we observe that Theorem 2.1 shows that the class of curves which separate $E$ and $F$ simply is non-empty. 

Now let $l$ be the infimum of the lengths of all geodesics which separate $E$ and $F$ simply. Similarly to the proof of Theorem 1.1, $l >0$ and we can take a sequence of closed geodesics $\gamma_n$ which separate $E$ and $F$ simply and whose lengths tend to $l$. Again we can find a subsequence $\gamma_{n_k}$ which converges uniformly to a geodesic $\gamma$ of length $l$ and as before we must have that $\gamma_{n_k} = \gamma$ for $k$ large enough, whence $\gamma$ must also separate $E$ and $F$ simply. 

Let us say that $E$ is the set contained in a single component of $\cbar \setminus \gamma$ and let us denote this component by $D$. Since $\gamma$ is connected, 
by Theorem 14.4 on Page 124 of \cite{New}, $\partial D$ is connected and so 
$D$ must be simply connected. By Lemma 2.2, $\partial D$ consists of finitely many segments of hyperbolic geodesics. $\partial D$ is thus locally connected and an inverse Riemann map $\phi$ from $\D$ to $D$ will extend continuously to the boundary (e.g.  \cite{Mil} Theorem 17.14). 

Since by Lemma 2.2 $\gamma$ may only intersect itself at finitely many points, it is not hard to see that no segment of $\gamma$ can be bounded on both sides by $D$ whence $\partial D$ is no longer than $\gamma$. Also , we can use $\phi$ to assign an orientation to the segments of $\partial D$ (which may require that we reverse the orientation of some of these original segments of $\gamma$). Finally, by approximating $\partial D$ with the image under $\phi$ of a circle ${\mathrm C}(0,r)$ whose radius $r$ is close to $1$ and using the continuous extension of the Riemann map to the boundary of $\D$ to give a homotopy in $U$, we can conclude that $n(\partial D, z) = 1$ for points $z \in D$ and 
$n(\partial D, z) = 0$ for points $z$ outside $\overline D$. 

If our original curve $\gamma$ was simple, then clearly, $\partial D = \gamma$. Any meridian in the homology class of $\gamma$ must also separate $E$ and $F$ simply. 
By the minimality of the length of $\gamma$, it follows from Theorem 1.3 that $\gamma$ must be a meridian. 

Now suppose $\gamma$ contains self-intersections, let $z$ be a point on $\partial D$ which is not a point of self-intersection of $\gamma$ and consider $\partial D$ as a curve beginning and ending at $z$. By Theorem 2.2, we can find a curve 
$\alpha$ in the fixed endpoint homotopy class of $\partial D$ which is a closed geodesic arc and which is smooth except possibly for a `corner' at $z$. Since $\gamma$ contains self-intersections, some of which must lie on $\partial D$, $\alpha$ must be shorter than $\partial D$. If we now apply Theorem 1.3 to obtain a meridian $\tilde \gamma$ in the homology class of $\alpha$ and thus of $\partial D$, $\tilde \gamma$ must be shorter than $\partial D$ and thus than $\gamma$ which would contradict the assumption that $\gamma$ was as short as possible. $\Box$
 
We remark that is a classical fact that provided a hyperbolic domain has a systole, i.e. a shortest geodesic, then, unless the the domain is a thrice punctured sphere, it must be a simple curve. The proof is a simple consequence of Theorem 4.2.4 on Page 100 of \cite{Bus} and the above. 

We now turn to proving the results we will need for Theorem 1.5.

\begin{theorem} (Nesting Theorem for Simple Closed Geodesics) Let $\gamma_1$, $\gamma_2$ be two simple closed geodesics in a hyperbolic domain $U$ and suppose they separate $\cbar \setminus U$ into sets $E_1$, $F_1$ and $E_2$, $F_2$ respectively. Suppose further that $E_1$ is connected and strictly contained in $E_2$, so that $\gamma_1$ is a principal meridian. Then $\gamma_1$ is contained in the component of $\cbar \setminus \gamma_2$ which contains $E_2$. In particular, $\gamma_1$ and $\gamma_2$ are disjoint. 
\end{theorem}

\proof  Let $U_1$, $V_1$ and $U_2$, $V_2$ be the complementary components of $\gamma_1$, $\gamma_2$ respectively which we label so that $E_1 \subset U_1$, $E_2 \subset U_2$. 

Since $E_1, E_2 \setminus E_1$ are both in $U_2$, $U_2$ contains points in $U_1$ and points in $V_1$ and so $\gamma_1$ must meet $U_2$. We thus need to show that $\gamma_1 \subset U_2$.  

So suppose now this fails and $\gamma_1 \cap (\cbar \setminus U_2) \neq \emptyset$. Since there is only one geodesic through a given point in a given direction, $\gamma_1$ and $\gamma_2$ cannot be tangent to each other at a point of intersection. Hence we can assume that there are points of $\gamma_1$ in $V_2$ and also that $U_1$ meets $\gamma_2$ and $V_2$.

Since $\gamma_2$ meets $U_1$, we can find an arc $\alpha_2$ of $\gamma_2$ which gives a crosscut of $U_1$ and by Theorem 11.7 on Page 118 of \cite{New}, this divides $U_1$ into two subdomains which we will call $A_1$ and $B_1$ as well as dividing $\gamma_1$ into two arcs. Each of these arcs meets one and only one of these subdomains since otherwise one of them would have a boundary which is contained in $\alpha_2$ which is clearly impossible as $\alpha_2$ is an arc with distinct endpoints and $\cbar \setminus \alpha_2$ is thus an unbounded simply connected domain. Since $E_1$ is connected, it must lie in just one of these domains, say, $A_1$. 

$B_1$ is then a Jordan domain which is entirely contained in $U$ with $\partial B_1$ consisting of $\alpha_2$ and an arc $\alpha_1$ of $\gamma_1$. It is then routine to find a slightly larger simply connected domain $\tilde B_1 \subset U$ such that $\alpha_1$ and $\alpha_2$ are fixed endpoint homotopic in $\tilde B_1$ and thus in $U$. But this is impossible by Theorem 2.2 as there is only one geodesic in any given homotopy class joining any two specified points of $U$. With this contradiction, the proof is complete. $\Box$

We remark here that once we know the meridians $\gamma_1$ and $\gamma_2$ are disjoint, we know that by the collar lemma we can find disjoint collars about each curve whose thicknesses in the hyperbolic metric of $U$ are uniformly bounded in terms of the lengths of these curves. 

\begin{lemma} Let $U$ be hyperbolic domain and let $\gamma_1$, $\gamma_2$ be two simple closed geodesics in $U$ which are homologous in $U$ and suppose that one of these curves lies in the closure of one of the complementary components of the other. Then $\gamma_1 = \gamma_2$. 
\end{lemma}

\proof As before, let $U_1$, $U_2$ and $V_1$, $V_2$ denote the bounded and unbounded complementary components of $\gamma_1$ and $\gamma_2$ respectively.
Since $\gamma_1 \underset{U} \approx \gamma_2$ and these curves are both simple, they separate $\cbar \setminus U$ into the same disjoint closed subsets $E$ and $F$. Using our earlier conventions, we may assume without loss of generality that $\infty \in F$ and that $n(\gamma_1, z) = n(\gamma_2, z) =1$ for all $z \in E$ while $n(\gamma_1, z) = n(\gamma_2, z) =0$ for all $z \in F$. Further, in order to be definite, we may assume without loss of generality that $\gamma_1$ lies inside $\overline U_2$, the proof in the other cases being essentially identical. If these curves are different, then, as before, any intersections must be transversal. Thus we can assume that either they are the same or $\gamma_1$ lies inside $U_2$.

So suppose these curves are different in which case they are disjoint. By the collar lemma (\cite{KL} Page 148, Lemma 7.7.1) we can find conformal annuli $A_1$ and $A_2$ whose equators are $\gamma_1$ and $\gamma_2$ respectively and which we can make sufficiently narrow about these curves so that they are disjoint. If we then  
let $\tilde A = \tilde A_2  \cup (U_2 \setminus U_1) \cup \tilde A_1$, $\tilde A$ is a conformal annulus whose boundary consists of two disjoint smooth simple closed curves. 

Let $a$ be a point in $E$ so that $n(\gamma_1, a)= n(\gamma_2, a) = 1$, let $V$ be the inverse image of $U$ under the covering map $e^z + a$ and let $\tilde S$ be the inverse image of $\tilde A$ under the same map. $\tilde S$ is then the region bounded between two vertically infinite simple curves which are invariant under a vertical shift by $2\pi i$ and in particular is simply connected. $\gamma_1$ and  $\gamma_2$ lift under $e^z + a$ to curves $\eta_1$ and $\eta_2$ respectively. These curves are also invariant under a vertical shift by $2\pi i$ and going around each of $\gamma_1$, $\gamma_2$ corresponds to the same covering transformation $z \mapsto z + 2\pi i$. 

Since $U$ is hyperbolic and $V$ is a covering of $U$, it follows that $V$ must also be hyperbolic and we can  
let $\psi : \D \mapsto V$ be a universal covering map. Since $\tilde S$ is simply connected any locally defined branch of $\psi^{\circ -1}$ can be defined on all of $\tilde S$ by the monodromy theorem. Thus $\eta_1$ and $\eta_2$ can be lifted to two geodesics in $\D$ which by Cauchy's theorem are a bounded hyperbolic distance apart and going around each of $\gamma_1$ and $\gamma_2$ can be made to correspond to the the same lifting under $\psi$ of 
$z \mapsto z + 2\pi i$. However, since these two geodesics in $\D$ must then be the same, it follows that $\gamma_1 = \gamma_2$ and, with this contradiction, the result follows. $\Box$

{\bf Proof of Theorem 1.5\hspace{.4cm}} Without loss of generality we can assume that $E$ is connected. Suppose now that the first part of the statement is false and we can find two different simple closed geodesics $\gamma_1$ and $\gamma_2$ which separate $E$ and $F$ and which we can then assume are homologous in $U$. Any intersection of these two curves must be transversal and if we let $U_1$ and $V_1$ be the bounded and unbounded components of the complement of $\gamma_1$, then $\gamma_2$ must meet both $U_1$ and $V_1$.  
Essentially the same argument as used in the proof of Theorem 2.7 now shows that 
this is impossible and so $\gamma_1$ and $\gamma_2$ are disjoint. But this is impossible by Lemma 2.3 above and, with this contradiction, we see that there can only be one such curve. 

For the second part of the statement, since at least one of $E$, $F$ is connected, any curve which separates them must do so simply. Suppose now that $\gamma$ is the meridian which separates $E$ and $F$. If $\tilde \gamma$ is another geodesic which separates these sets, then we have already seen that $\tilde \gamma$ cannot be simple and since this curve must separate $E$ and $F$ simply, by Theorem 2.5, it must be longer than $\gamma$. The result then follows. $\Box$

\section{Finitely Connected Domains}

Let $U$ be a hyperbolic domain of connectivity $n$ and let $K^1, K^2, \ldots \ldots, K^n$ denote the complementary components of $U$ where we will always assume that $K^n$ is the unbounded one (Ahlfors uses the same convention in \cite{Ahl}). 

One can then see using elementary combinatorics that up to homology there are at most $E(n) := 2^{n-1} - 1$ classes of meridians for $U$. Note that this includes the simply connected case where there is just one complementary component ($U = \C$ not being allowed as we only consider hyperbolic domains) in which case $n =1$ and we have no meridians. 

We say $U$ is \emph{non-degenerate} if none of the complementary components are points. If $n \ge 2$ so that $E(n) > 0$, then by Theorem 2.3, it is easy to see that we have exactly $E(n)$ meridians if and only if $U$ is non-degenerate. It is also easy to check that the maximum possible number of principal meridians is given by $P(n):= \min\{n,E(n)\}$. 
Moreover, those meridians which fail to exist when some of the complementary components are points are principal meridians and so there must in fact be at least $E(n) - P(n)$ meridians. 

By Theorem 1.4, the principal meridians of $U$ are uniquely defined. If we can find $P(n)$ principal meridians, let us call such a collection the \emph{principal system of meridians} or simply the \emph{principal system} for $U$. If we can find a full collection of $E(n)$ meridians, let us call such a collection an \emph{extended system of meridians} or simply an \emph{extended system} for $U$. We have shown the following.

\begin{proposition}
If $U$ is a domain of finite connectivity $n \ge 2$, then $U$ has at least $E(n) - P(n)$ meridians and any principal meridians of $U$ which exist are uniquely defined. Furthermore, the following are equivalent:
\vspace{0cm}

\begin{enumerate}
\item $U$ is non-degenerate;

\vspace{.25cm}
\item $U$ has $P(n)$ principal meridians;

\vspace{.25cm}
\item $U$ has $E(n)$ meridians in distinct homology classes.

\end{enumerate}

\end{proposition}

Note that if $n =2$ or $3$, than any meridians of $U$ which exist must be principal. The first case where we can have meridians which are not principal is $n=4$ as can be seen in Figure 1. If none of the complementary components except possibly the unbounded one are points, then the $n-1$ principal meridians which separate the bounded complementary components from the rest of $\cbar \setminus U$ give us a homology basis for $U$ (note that if $\infty \in U$, then we will need all $P(n)$ principal meridians for our homology basis). This basis is efficient in the sense that it contains the fewest possible curves and that these curves are as short as possible in their homology classes. As an added advantage, by Theorem 2.7, these principal meridians are disjoint both from each other and from any other meridians of the domain. 

Finally, we remark that the principal meridians are useful for estimating the separation between the components of $\cbar \setminus U$. This can be used to give a measure of how close $U$ comes to being degenerate. This relates to questions regarding the \Car topology for convergence of pointed domains, in particular whether or not a limit of $n$-connected pointed domains is again $n$-connected and we propose to examine this situation in detail in subsequent papers.


\end{document}